# GREAT CIRCLE FIBRATIONS and CONTACT STRUCTURES on ODD-DIMENSIONAL SPHERES

**Herman Gluck and Jingye Yang**

**Abstract.** *It is known that for every smooth great circle fibration of the 3-sphere, the distribution of tangent 2-planes orthogonal to the fibres is a contact structure, in fact a tight one, but we show here that, beginning with the 5-sphere, there exist smooth great circle fibrations of all odd-dimensional spheres for which the tangent hyperplane distribution orthogonal to the fibres is* ***not*** *a contact structure.*

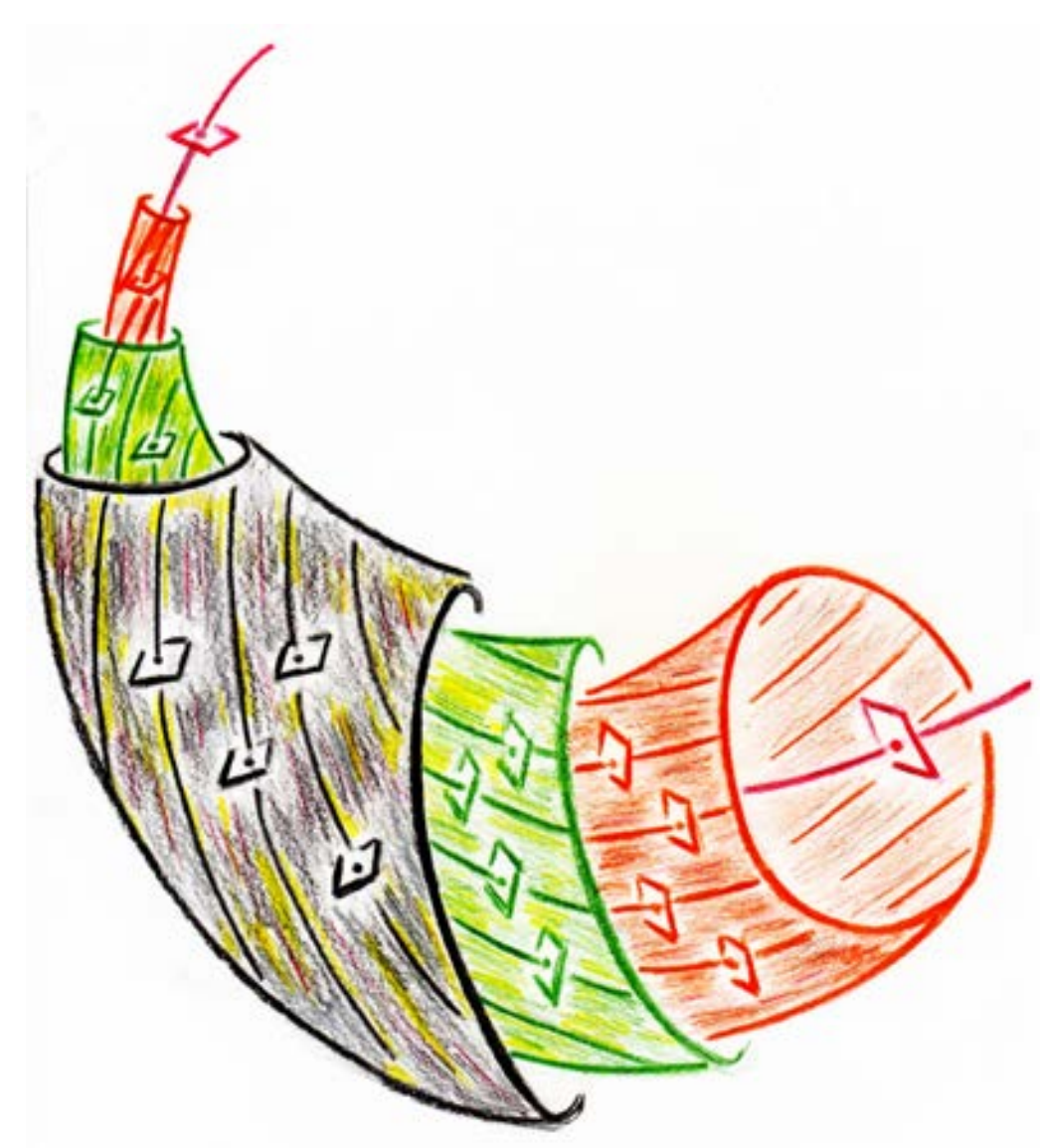

**THEOREM**. ***On all odd-dimensional spheres, beginning with the 5-sphere, there exist smooth great circle fibrations for which the tangent hyperplane distribution orthogonal to the fibres is not a contact structure.***

***By contrast, on the 3-sphere, for every smooth great circle fibration, the tangent hyperplane distribution orthogonal to the fibres is a tight contact structure.***

# INTRODUCTION

## Key propositions, informally phrased.

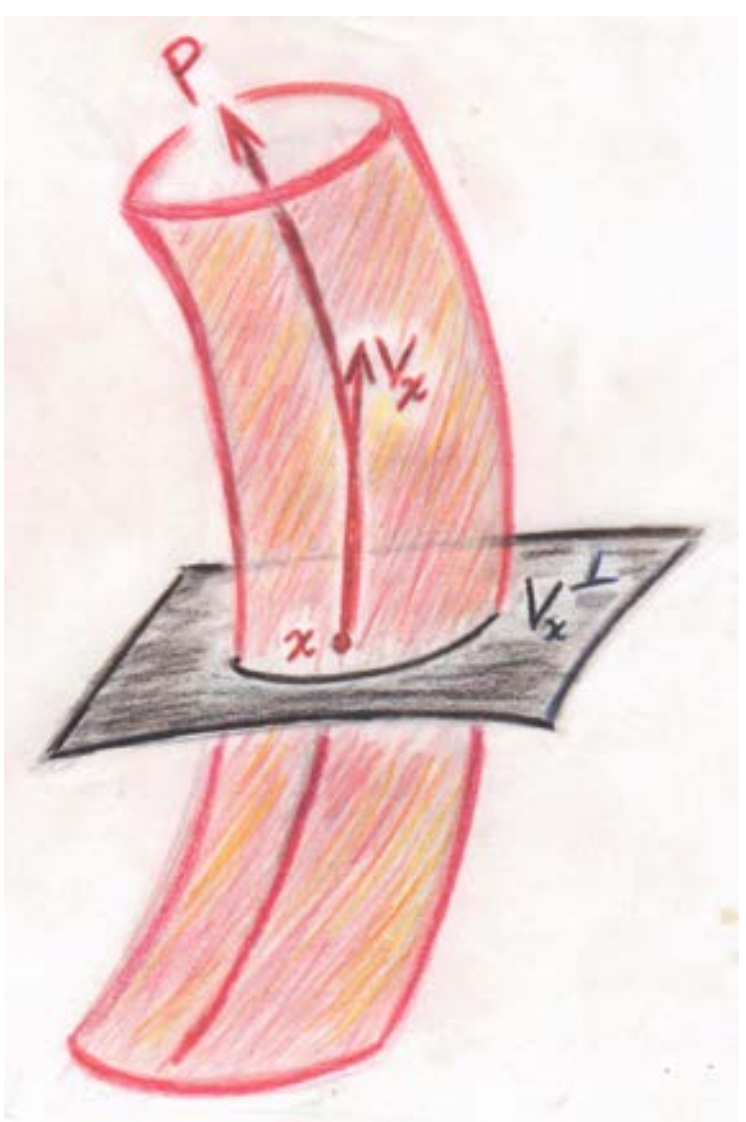

The unit tangent vector $V_x$ at the point $x$ on a fibre $P$ of a smooth fibration of $S^{2n+1}$ by oriented great circles, and the tangent 2n-plane $V_x^{\perp}$ orthogonal there to $V_x$ .

**PROPOSITION 1.** ***Given a smooth fibration F of $S^{2n+1}$ by oriented great circles with unit vector field V tangent to the fibres , there is at each point $x \in S^{2n+1}$ a linear map $T_x : V_x^{\perp} \rightarrow V_x^{\perp}$ which expresses "twisting" of the fibres about the one through $x$ , and which has no real eigenvalues.***

**PROPOSITION 2.** ***A necessary and sufficient condition for the orthogonal tangent 2n-plane distribution $\xi_F = \{V_x^{\perp} : x \in S^{2n+1}\}$ to be a contact structure is that the skew-symmetrization $T_x - T_x^{tr}$ of this twisting map be nonsingular at each point $x \in S^{2n+1}$ .***

This condition is always satisfied on the 3-sphere, but can be violated beginning on the 5-sphere.

## Picturing everything inside a Grassmann manifold.

Given a fibration $F$ of $S^{2n+1}$ by oriented great circles, each fibre $P$ of $F$ lies in and orients some 2-plane through the origin in $R^{2n+2}$. We will denote this 2-plane by $P$ as well, and view it as a single point in the Grassmann manifold $G_2R^{2n+2}$ of all such oriented 2-planes. The base space $M_F$ of $F$ then appears as a 2n-dimensional topological submanifold of $G_2R^{2n+2}$, and if the fibration $F$ is smooth, then the submanifold $M_F$ is also smooth.

Given an oriented 2-plane $P$ through the origin in $R^{2n+2}$, let $P^\perp$ denote the 2n-plane through the origin orthogonal to it. There is no need to orient $P^\perp$.

The 4n-dimensional vector space $Hom(P, P^\perp)$ serves simultaneously as a large coordinate neighborhood about $P$ in $G_2R^{2n+2}$, and as the tangent space $T_P(G_2R^{2n+2})$ to this Grassmann manifold, as follows.

Suppose that the oriented 2-plane $P'$ in $R^{2n+2}$ contains no vector orthogonal to $P$, and suppose that its orthogonal projection to $P$ is orientation-preserving. Let $N(P)$ be the collection of all such 2-planes $P'$. This set $N(P)$ is the domain of a coordinate chart, as follows. Given $P'$ in $N(P)$, we view $P'$ as the graph of a linear transformation $L_{P'} : P \to P^\perp$, and we match $P'$ with $L_{P'}$. In this way, $P$ itself is matched with the zero transformation. Vice versa, if we start with a linear map $L: P \to P^\perp$, then its graph is the oriented 2-plane $P_L$.

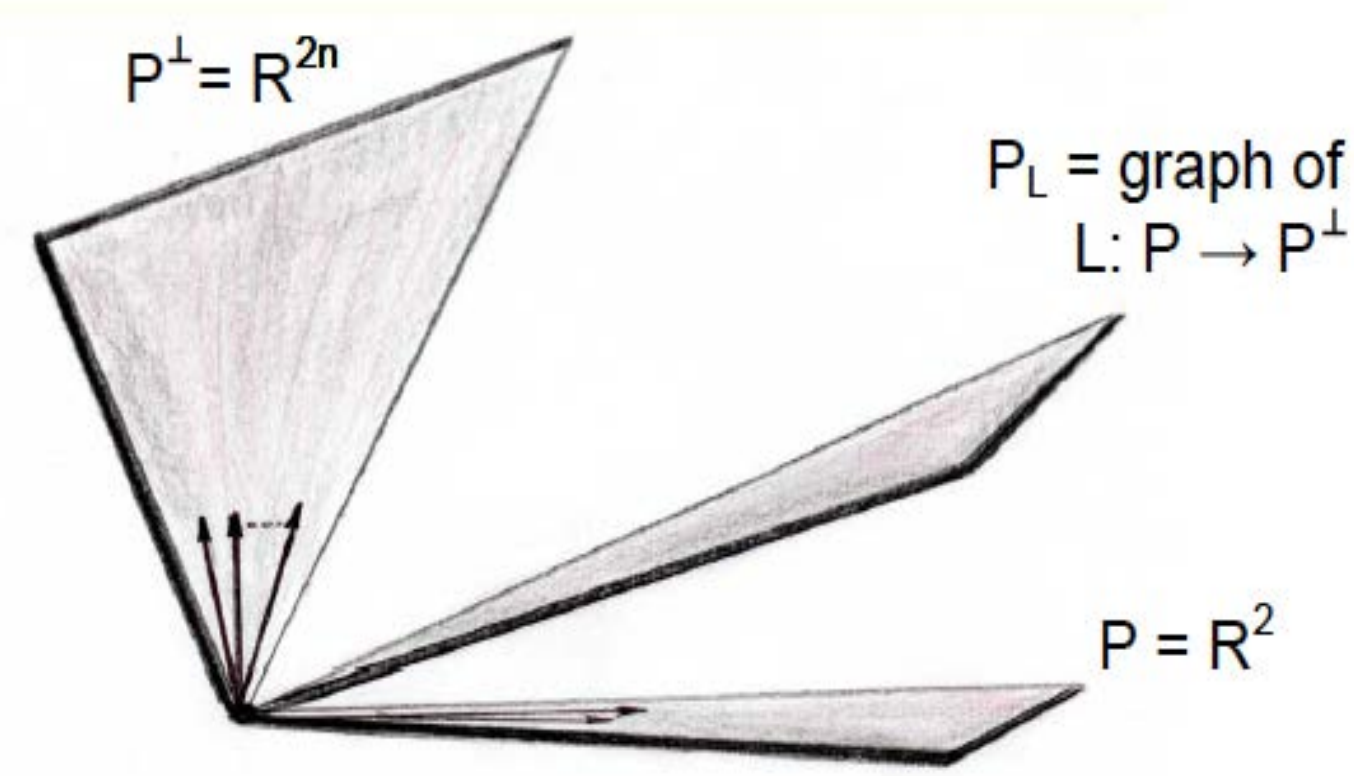

Let $e_1$, $e_2$ be an ordered orthonormal basis for $P$ consistent with its orientation. There is a circle's worth of these.

An element $L$ of $Hom(P, P^\perp)$ determines and is determined by the two vectors $L(e_1)$ and $L(e_2)$ in $P^\perp$, and hence $Hom(P, P^\perp)$ is isomorphic to $P^\perp + P^\perp$.

There is a circle's worth of such decompositions of $Hom(P, P^\perp)$, corresponding to the circle's worth of ordered orthonormal bases for $P$.

## New tools.

The argument [Gluck 2022] for our theorem on the 3-sphere used special features available only in this low dimension: quaternion multiplication, isometry of the Grassmann manifold $G_2R^4$ of oriented 2-planes through the origin in $R^4$ with the product space $S^2 \times S^2$, and the description of the infinite-dimensional moduli space of all great circle fibrations of the 3-sphere as two copies of the space of strictly distance-decreasing maps from $S^2$ to $S^2$.

To prove the negative result on higher dimensional spheres, we will first show how to free ourselves from the above special tools by finding new more general ones.

Consider an oriented great circle fibration F of $S^{2n+1}$ which contains a fixed great circle fibre P. If P' is another great circle on $S^{2n+1}$ which intersects P, then because the fibres of F are disjoint, the base space $M_F$ in $G_2R^{2n+2}$ cannot also pass through P'. This motivates the following definition.

The ***bad set*** $BS(P) \subset G_2R^{2n+2}$ consists of all oriented 2-planes through the origin in $R^{2n+2}$ which meet P in at least a line. If $M_F$ contains the great circle fibre P, then $M_F$ intersects the bad set BS(P) only at P and nowhere else.

The ***bad cone*** $BC(P) \subset T_P(G_2R^{2n+2})$ is the tangent cone to the bad set at P.

**TOOL 1.** ***A closed connected smooth 2n-dimensional submanifold of*** $\mathbf{G_2R^{2n+2}}$ ***is the base space of a fibration of*** $\mathbf{S^{2n+1}}$ ***by oriented great circles if and only if it is transverse to the bad cone at each of its points.***

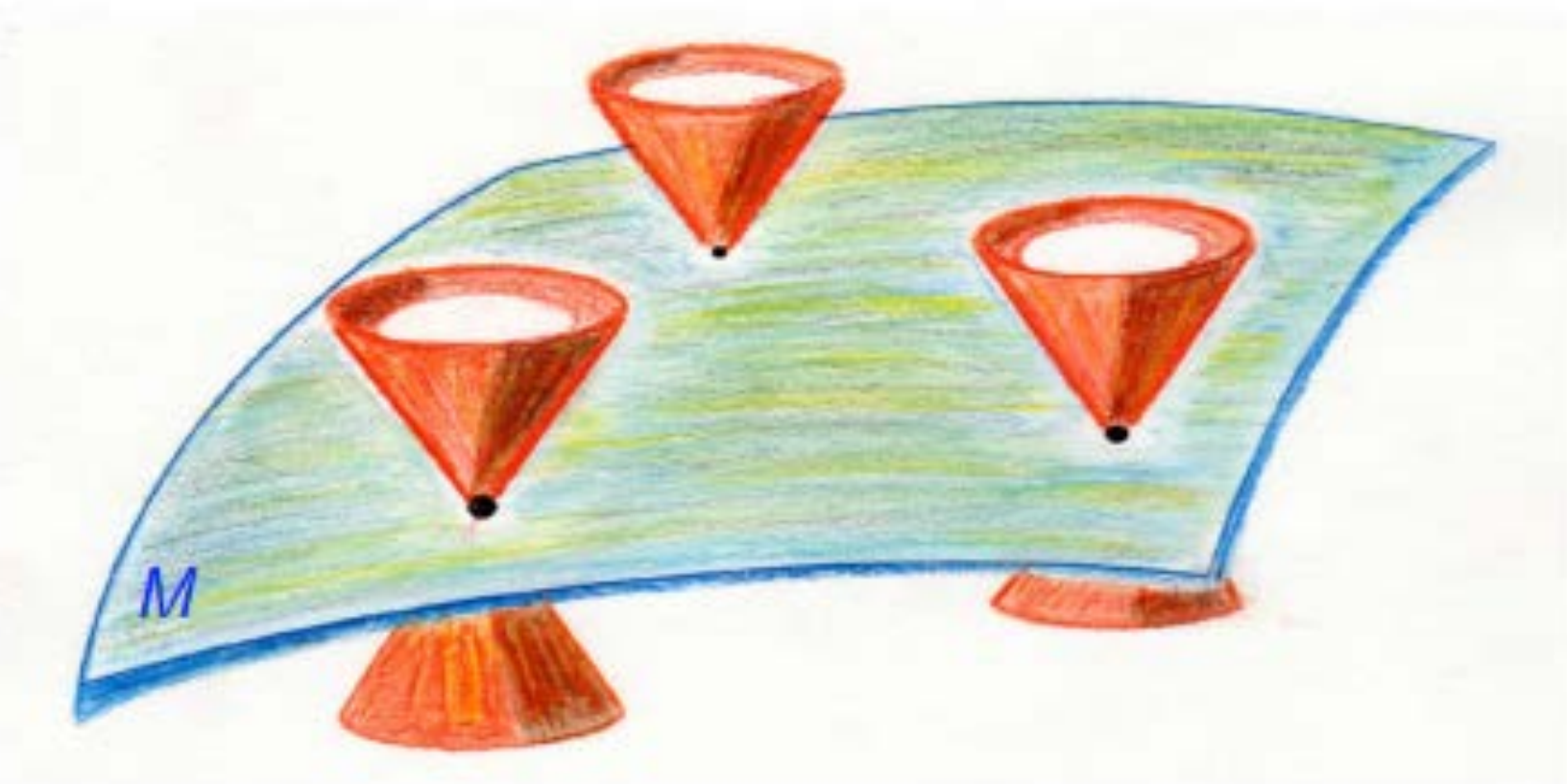

For $S^3$, this is Theorem B of [Gluck-Warner 1983]]. For smooth fibrations of spheres by great subspheres of any dimension, this is Theorem 4.1 of [Gluck-Warner-Yang 1983]. This was proved again for all great circle fibrations of $S^{2n+1}$ by Benjamin McKay [McKay 2004] from a different point of view.

**LOCAL VERSION OF TOOL 1.** ***A sufficiently small neighborhood of a point on a smooth 2n-cell in*** $G_2R^{2n+2}$ ***which is transverse to the field of bad cones is the base space of a fibration of a thin tube in*** $S^{2n+1}$ ***by oriented great circles [Cahn-Gluck-Nuchi 2018].***

**TOOL 2.** ***A 2n-dimensional linear subspace of*** **Hom(P, $P^\perp$)** $\cong$ **$P^\perp$ + $P^\perp$** ***is transverse to the bad cone*** **BC(P)** ***if and only if it is the graph of a linear map with no real eigenvalues from one*** **$P^\perp$** ***summand to the other.***

See [Cahn-Gluck-Nuchi 2018].

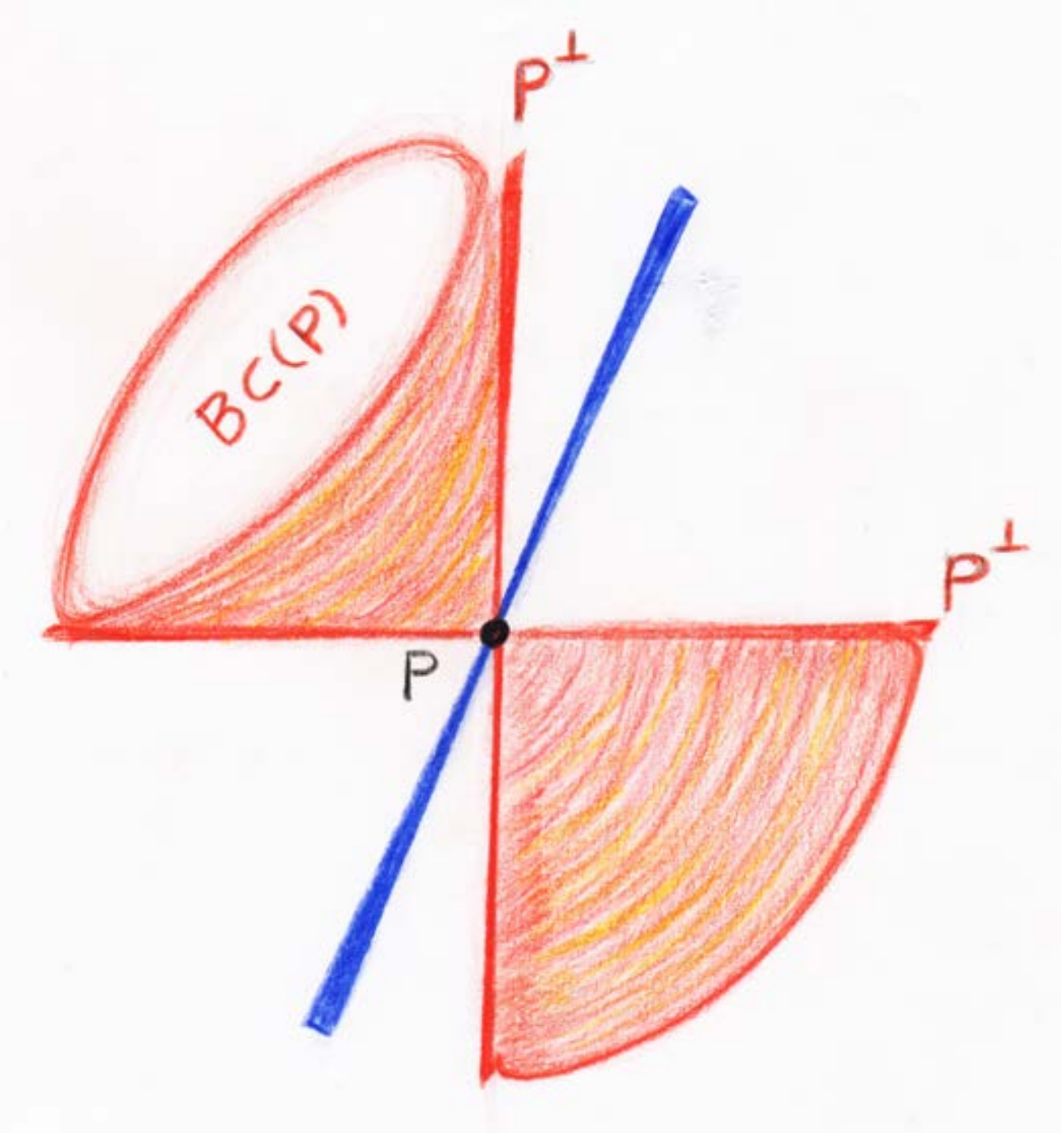

**A linear subspace transverse to the bad cone BC(P) in $T_P(G_2R^{2n+2})$**

A ***germ*** of a fibration of $S^{2n+1}$ by oriented great circles consists of such a fibration in an open neighborhood of a given fibre P , with two germs about P equivalent if they agree on some smaller neighborhood of P . To ***extend*** such a germ to a fibration of $S^{2n+1}$ means to find a fibration of $S^{2n+1}$ which agrees with the given germ on some neighborhood of P .

**TOOL 3.** ***Every germ of a smooth fibration of*** **$S^{2n+1}$** ***by oriented great circles extends to such a fibration of all of*** **$S^{2n+1}$.**

This is Theorem A of [Cahn-Gluck-Nuchi 2018].

## Remarks.

- The standard Hopf fibration $H$ of the 3-sphere is obtained by starting with an orthogonal complex structure $J$ on $C^2 = R^4$, and then intersecting the J-complex lines through the origin with the unit 3-sphere centered at the origin to obtain the great circle fibres. These great circles can be oriented in the direction of complex multiplication by $i$.

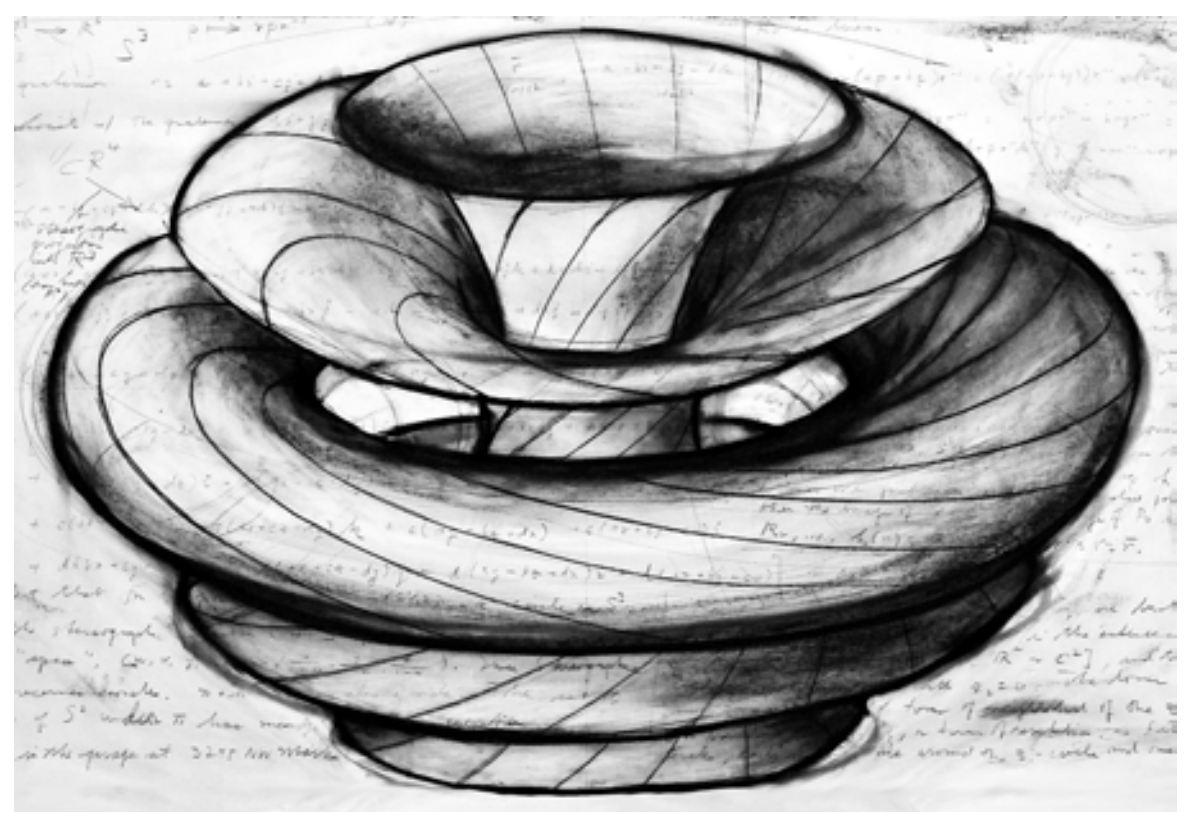

**Hopf fibration of the 3-sphere by great circles**
**Lun-Yi Tsai Charcoal and graphite on paper 2006**

The distribution $\xi_H$ of tangent 2-planes on $S^3$ which are orthogonal to these Hopf fibres is known as the ***standard tight contact structure*** on $S^3$.

- The subject of fibrations of round spheres by great subspheres, both in and of itself, and in the way it applies to the Blaschke Problem in Differential Geometry, is nicely summarized in [McKay 2004 and 2015].

- A overview of contact structures and contact geometry can be found in [Eliashberg 1992 and 1993], in [Etnyre 2003] and in [Geiges 2008].

## Acknowledgements.

- Many thanks to Jo Nelson, who gave a beautiful talk here at Penn in 2017, which inspired thinking about the possible existence of a connection between fibrations of spheres by great circles and contact geometry.

- Many thanks also to our colleagues Dennis DeTurck, Mona Merling, Leandro Lichtenfelz, and Ziqi Fang for helpful advice.

- Special thanks to the referee for suggesting coordinate-free, conceptual proofs of Propositions 1 and 2 ... much more satisfying than the ones we originally gave.

## GREAT CIRCLE FIBRATIgONS

### The twisting map.

Let $F : S^1 \subset S^{2n+1} - p \rightarrow M_F$ be a smooth fibration of $S^{2n+1}$ by oriented great circles, and let $V$ be the unit vector field on $S^{2n+1}$ tangent to the fibres of $F$ which points in the direction of their orientation. Fixing a point $x$ on $S^{2n+1}$, let $V_x$ denote the value of $V$ at $x$ and $V_x^{\perp}$ the tangent 2n-plane there orthogonal to $V_x$.

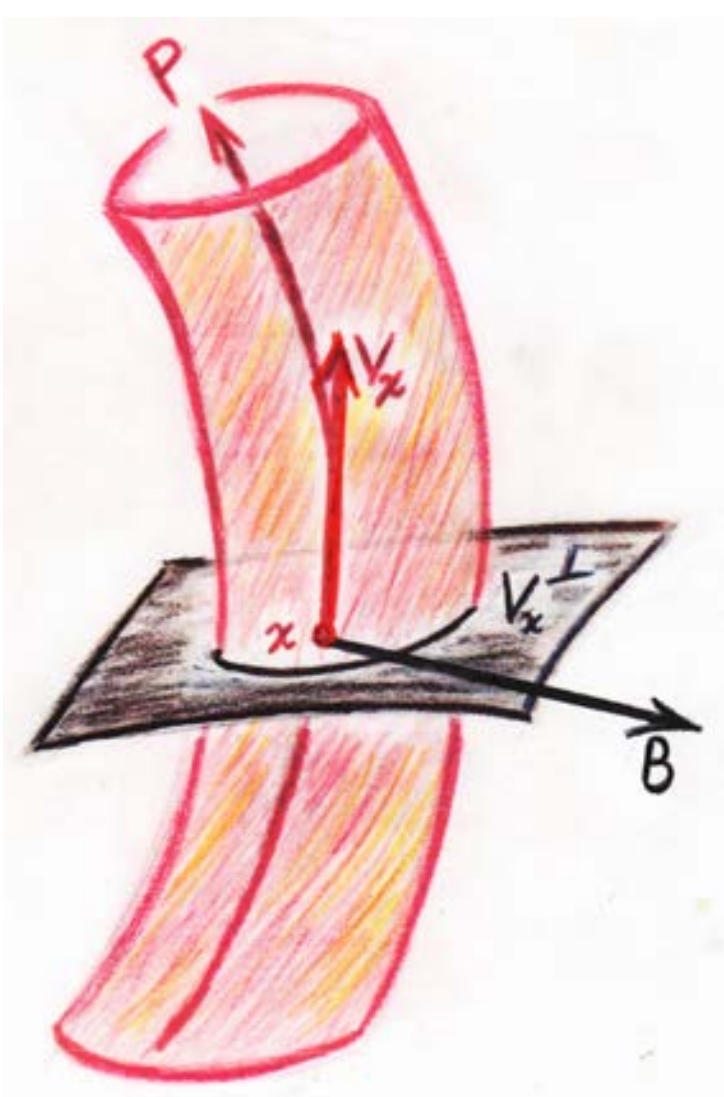

Keeping $x$ fixed, we will use $P$ to denote either the fibre of $F$ through $x$ or the oriented 2-plane through the origin which it spans. Parallel translation of $V_x^{\perp}$ to the origin takes it to the 2n-plane $P^{\perp}$ orthogonal to $P$.

Let $B$ be a unit vector in $V_x^{\perp}$, and consider the covariant derivative $\nabla_B V$, which is some tangent vector to $S^{2n+1}$ at $x$. See [do Carmo 1992, Chapter 2] for a discussion of covariant derivatives and their properties. But for simplicity, to compute the covariant derivative at $x$ of a vector field along a parametrized curve in $S^{2n+1}$, just compute its ordinary derivative there in $R^{2n+2}$ and project the answer orthogonally back tangent to $S^{2n+1}$ at $x$.

Note that $\nabla_B V$ also lies in $V_x^{\perp}$, because

$$\langle \nabla_B V , V \rangle = 1/2\, B \langle V , V \rangle = 0 ,$$

since $\langle V , V \rangle = 1$.

**The restriction $(\nabla V)|_{V_x^{\perp}}$, which takes $B \rightarrow \nabla_B V$, is a linear map $V_x^{\perp} \rightarrow V_x^{\perp}$ which measures the sidewise rate of change of the vector field $V$ at the point $x$.**

This is our ***twisting map*** $T_x = (\nabla V)|_{V_x^{\perp}}$.

## The twisting map for the Hopf fibrations.

Using the identification of $R^{2n+2}$ with $C^{n+1}$ to help us choose a nice basis for $P^\perp$, the twisting map $T : P^\perp \rightarrow P^\perp$ has the matrix representation

$$\begin{matrix} \mathbf{0} & \mathbf{1} & & & & & \\ \mathbf{-1} & \mathbf{0} & & & & & \\ & & \mathbf{0} & \mathbf{1} & & & \\ & & \mathbf{-1} & \mathbf{0} & & & \\ & & & & \cdots & & \\ & & & & & \mathbf{0} & \mathbf{1} \\ & & & & & \mathbf{-1} & \mathbf{0} \end{matrix} .$$

To confirm this for the Hopf fibration on $S^3$, regard this sphere as the space of unit quaternions, and on it consider the orthonormal basis of left-invariant vector fields:

$$A(x) = x\,i \;, \quad B(x) = x\,j \;, \quad C(x) = x\,k \,.$$

Computing Euclidean covariant derivatives in $R^4$ and projecting orthogonally back to $S^3$, we get:

$$\nabla_A B = C \,, \;\; \nabla_B C = A \,, \;\; \nabla_C A = B \quad \text{and} \quad \nabla_B A = -C \,, \;\; \nabla_C B = -A \,, \;\; \nabla_A C = -B \,.$$

Now choose the Hopf fibres so that the unit vector field $V$ along them lines up with $A$.

Then

$$\nabla_B V = \nabla_B A = -C \quad \text{and} \quad \nabla_C V = \nabla_C A = B \;,$$

which corresponds to the 2 x 2 matrix in the upper left corner above.

We leave confirmation for the Hopf fibrations on $S^{2n+1}$ to the reader.

## The Grassmann and Stiefel manifolds $G_2R^{2n+2}$ and $V_2R^{2n+2}$ .

Earlier we introduced the Grassmann manifold $G_2R^{2n+2}$ of oriented 2-planes through the origin in $R^{2n+2}$ and indicated its importance for us.

We will also use the Stiefel manifold $V_2R^{2n+2}$ of orthonormal 2-frames $(e_1 , e_2)$ in $R^{2n+2}$ . It is the total space of a bundle $S^1 \subset V_2R^{2n+2} \to G_2R^{2n+2}$ , with the projection map taking the orthonormal 2-frame $(e_1 , e_2)$ to the oriented 2-plane $[e_1 , e_2]$ through the origin which it spans.

A fibre of the projection map from $V_2R^{2n+2} \to G_2R^{2n+2}$ is the circle's worth of orthonormal 2-frames which span the same oriented 2-plane through the origin in $R^{2n+2}$ . If $(e_1 , e_2)$ is one point (orthonormal 2-frame) on a given fibre of the Stiefel bundle, then all the points on this fibre are $( e_1 \cos t + e_2 \sin t , - e_1 \sin t + e_2 \cos t )$ , $0 \le t \le 2\pi$ . In particular, the vector $W = (e_2 , - e_1)$ is tangent to this Stiefel fibre at the point $(e_1 , e_2)$ .

The Stiefel manifold is important to us because we will compute velocities along curves in the Grassmann manifold by first lifting them to the Stiefel manifold, computing velocities there, and then projecting the result back down to the Grassmann manifold.

## Tangent spaces to these manifolds.

Let $P$ be an oriented 2-plane through the origin in $R^{2n+2}$ and $P^\perp$ its orthogonal complement. We saw earlier how $Hom(P, P^\perp)$ can serve as the tangent space $T_P G_2R^{2n+2}$ to our Grassmann manifold at $P$ , and that if we choose an ordered orthonormal basis $e_1 , e_2$ for $P$ , then linear maps $L \in Hom(P, P^\perp)$ are in one-to-one correspondence with pairs of vectors $(L(e_1) , L(e_2))$ in $P^\perp + P^\perp$ .

There is a corresponding decomposition of the tangent space to the Stiefel manifold $V_2R^{2n+2}$ at the point $(e_1 , e_2)$ . The tangent vector $(e_2 , - e_1)$ spans the direction of the Stiefel fibre. Independent of this, $e_1$ can tend to move in the directions of $P^\perp$ while $e_2$ stays fixed, and alternatively, $e_2$ can tend to move in the directions of $P^\perp$ while $e_1$ stays fixed.

The projection map $V_2 R^{2n+2} \to G_2 R^{2n+2}$ kills the direction of the Stiefel fibre, but its differential matches the $P^\perp + P^\perp$ tangent decomposition upstairs in the Stiefel manifold with that downstairs in the Grassmann manifold.

Using the Riemannian metrics from the next section, the projection map $V_2R^{2n+2} \to G_2R^{2n+2}$ is a Riemannian submersion, meaning that it preserves lengths orthogonal to the circle fibres.

## Riemannian metrics on the Stiefel and Grassmann manifolds.

We view $V_2R^{2n+2} \subset S^{2n+1} \times S^{2n+1}$ , and give the Stiefel manifold the Riemannian metric induced via this inclusion from the usual product of round metrics on the factor spheres.

Then the diagonal action of the special orthogonal group SO(2n+2) on $S^{2n+1} \times S^{2n+1}$ takes $V_2R^{2n+2}$ to itself via isometries, acts transitively on this Stiefel manifold except when $n = 0$ , and lets us write $V_2R^{2n+2} = SO(2n+2) / SO(2n)$ .

The standard Riemannian metric on the Grassmann manifold $G_2R^{2n+2}$ was explored in [Leichtweiss 1961] and [Wong 1967], see also [Gluck-Warner 1983]. The distance $d(P, Q)$ between two oriented 2-planes P and Q in $R^{2n+2}$ is defined in terms of the principal angles $\theta_1(P, Q)$ and $\theta_2(P, Q)$ between them by the formula

$$d(P, Q)^2 = \theta_1(P, Q)^2 + \theta_2(P, Q)^2 .$$

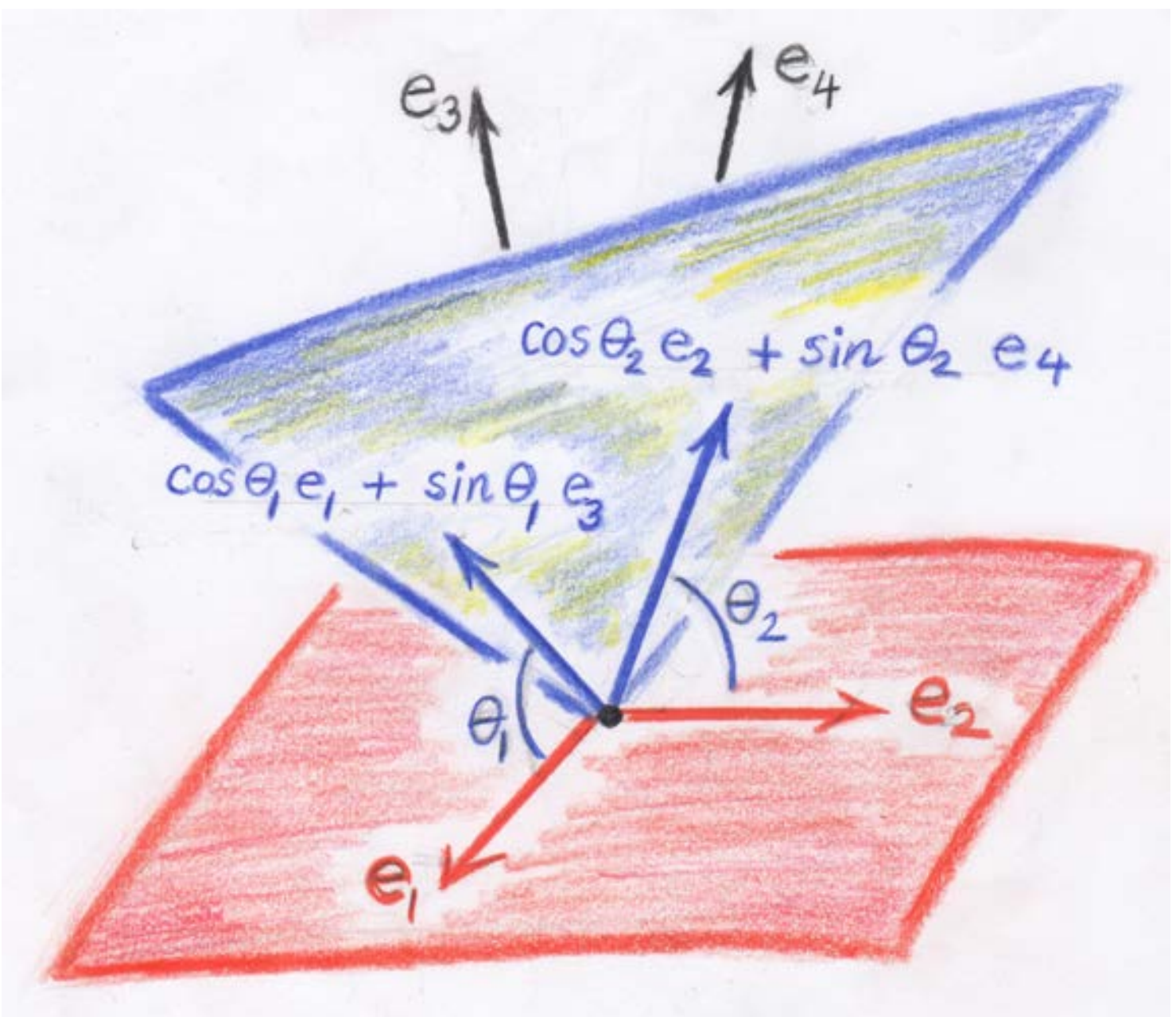

In the 4-space spanned by this figure, steady rotations of the $e_1 e_3$ and $e_2 e_4$ -planes within themselves at speeds proportional to the principal angles $\theta_1$ and $\theta_2$ define a geodesic in SO(4) which will move the 2-plane P to Q along a geodesic in $G_2R^4$ .

This action is covered by a horizontal geodesic in $V_2R^4$ which takes the orthonormal 2-frame $(e_1 , e_2)$ to the orthonormal 2-frame $(\cos\theta_1 e_1 + \sin\theta_1 e_3 , \cos\theta_2 e_2 + \sin\theta_2 e_4 )$ . Here "horizontal" means orthogonal to the fibres of the Stiefel projection, and this is confirmed at time zero, because a tangent vector to the Stiefel fibre through $(e_1 , e_2)$ at time zero is $(e_2 , -e_1)$ , while a tangent vector to the curve in question is $(\theta_1 e_3 , \theta_2 e_4)$ .

The above is representative of geodesics in $G_2R^{2n+2}$ and $V_2R^{2n+2}$ , since in both cases the geodesic connecting two elements of either of these spaces runs within a $G_2R^4$ or $V_2R^4$ .

The special orthogonal group SO(2n+2) acts isometrically and transitively on $G_2R^{2n+2}$ , and lets us write $G_2R^{2n+2} = SO(2n+2) / [SO(2n) \times SO(2)]$ .

**PROPOSITION 1. Given the smooth fibration $F : S^1 \subset S^{2n+1} \rightarrow M_F$ of the unit (2n+1)-sphere $S^{2n+1}$ by oriented great circles, the graph of the twisting map $T_x = (\nabla V)|_{V_x^{\perp}} : V_x^{\perp} = P^{\perp} \rightarrow V_x^{\perp} = P^{\perp}$ is the tangent 2n-plane $T_P M_F$ to the base space $M_F \subset G_2R^{2n+2}$ at the point $P$, and therefore $T_x$ can have no real eigenvalues.**

**Proof.**

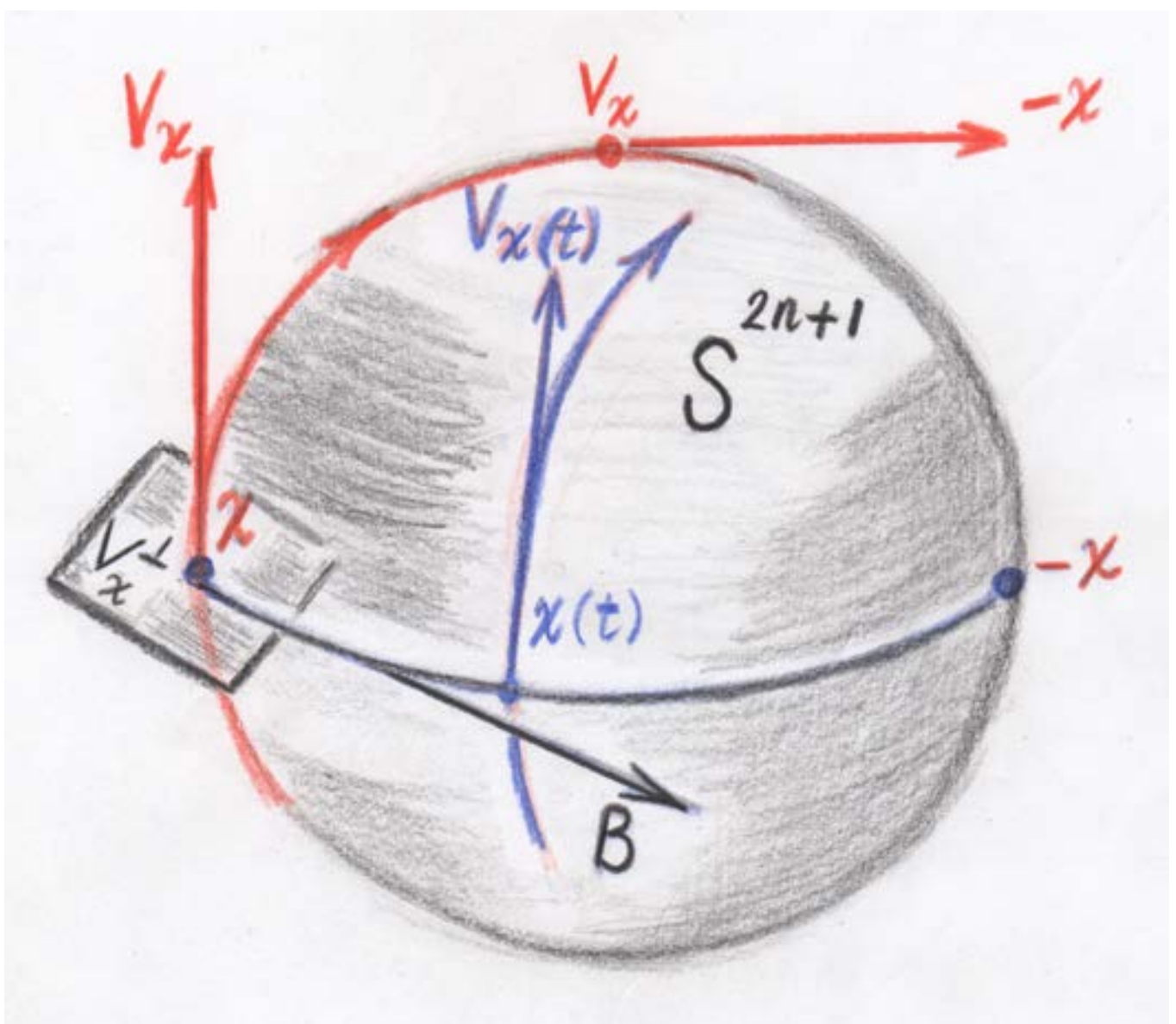

At the left in the figure above, we show the point $x$ in $S^{2n+1}$ and the unit vector $V_x$, tangent there to the oriented great circle fibre of $F$ through $x$, and also show the tangent 2n-plane $V_x^{\perp}$ orthogonal there to $V_x$. This great circle fibre of $F$ spans the oriented 2-plane $P$ through the origin in $R^{2n+2}$, while the tangent 2n-plane $V_x^{\perp}$, when translated to the origin, becomes the 2n-plane $P^{\perp}$.

If we parallel translate the unit vector $V_x$ at $x$ so that it begins at the origin, then it will end at the point $V_x$ shown at the top of the figure, still on the same fibre of $F$.

The tangent vector to this fibre at the point $V_x$ can be labeled $-x$, because when parallel translated to begin at the origin, it will end at the point $-x$, also on the same fibre of $F$.

The pair of vectors $(x, V_x)$ is an orthonormal 2-frame in $R^{2n+2}$, hence a point of the Stiefel manifold $V_2R^{2n+2}$.

It projects down to the oriented 2-plane $P = [x, V_x]$ in the Grassmann manifold $G_2R^{2n+2}$.

In $V_x^\perp$ we now pick a random unit vector $B$ and, as shown in the figure, move $x$ to the right along the great circle through $x$ and $B$, so that at time $t$ it is at the location $x(t) = x \cos t + B \sin t$. At that point, we have the unit vector $V(x(t))$ tangent to the fibre of $F$ through $x(t)$.

This gives us a smooth curve $\bar{\gamma}(t) = ( x(t) , V(x(t)) )$ in the Stiefel manifold $V_2R^{2n+2}$ which covers the smooth curve $\gamma(t) = [ x(t) , V(x(t)) ]$ in the Grassmann manifold $G_2R^{2n+2}$.

The initial velocity of the curve $\bar{\gamma}(t) = ( x(t) , V(x(t)) )$ in the Stiefel manifold is the pair of vectors $(d/dt)|_{t=0}\, x(t)$ and $(d/dt)|_{t=0}\, V(x(t))$.

By construction, $(d/dt)|_{t=0}\, x(t) = B$, which lies in $V_x^\perp$ and is tangent to $S^{2n+1}$ at $x$.

We claim that $(d/dt)|_{t=0}\, V(x(t)) = \nabla_B V$, with no projection to $S^{2n+1}$ needed.

By definition of $\nabla_B V$, it is the tangential component of $(d/dt)|_{t=0}\, V(x(t))$, so we need to show that there is no normal component of this vector, namely that the inner product $< x , (d/dt)|_{t=0}\, V(x(t)) > = 0$ at the point $x \in S^{2n+1}$.

To do this, we start with the formula $< x(t) , V(x(t)) > = 0$, since $V(x(t))$ is tangent to $S^{2n+1}$ at the point $x(t)$, and hence orthogonal to the vector $x(t)$.

We differentiate this formula with respect to $t$ at $t = 0$:

$$(d/dt)|_{t=0} < x(t) , V(x(t)) > = 0.$$

Applying the Leibniz Rule, we get

$$(d/dt)|_{t=0} < x(t) , V(x(t)) > = < (d/dt)|_{t=0}\, x(t) , V_x > + < x , (d/dt)|_{t=0}\, V(x(t)) > = 0 .$$

Now $(d/dt)|_{t=0}\, x(t) = B$ by construction, and we chose $B$ so that $< B , V_x > = 0$.

Hence $< x , (d/dt)|_{t=0}\, V(x(t)) > = 0$, as desired, and therefore

$$(d/dt)|_{t=0}\, V(x(t)) = \nabla_B V ,$$

as claimed.

So the initial velocity of the curve $\bar{\gamma}(t) = ( x(t) , V(x(t)) )$ in the Stiefel manifold $V_2R^{2n+2}$ is the pair of vectors

$$( (d/dt)|_{t=0}\, x(t) , (d/dt)|_{t=0}\, V(x(t)) ) = ( B , \nabla_B V ) .$$

We claim this velocity is orthogonal to the Stiefel fibre through $( x , V_x )$ .

We recall that the vector $( V_x , -x )$ is tangent to this Stiefel fibre at the point $( x , V_x )$ .

But then the initial velocity vector of the curve $\bar{\gamma}(t) = ( x(t) , V(x(t)) )$ at this point in the Stiefel manifold, namely the pair of vectors $( B , \nabla_B V )$ , is orthogonal at the point $( x , V_x )$ to the Stiefel fibre because

$$< B , V_x > + < \nabla_B V , -x > = 0 + 0 = 0 .$$

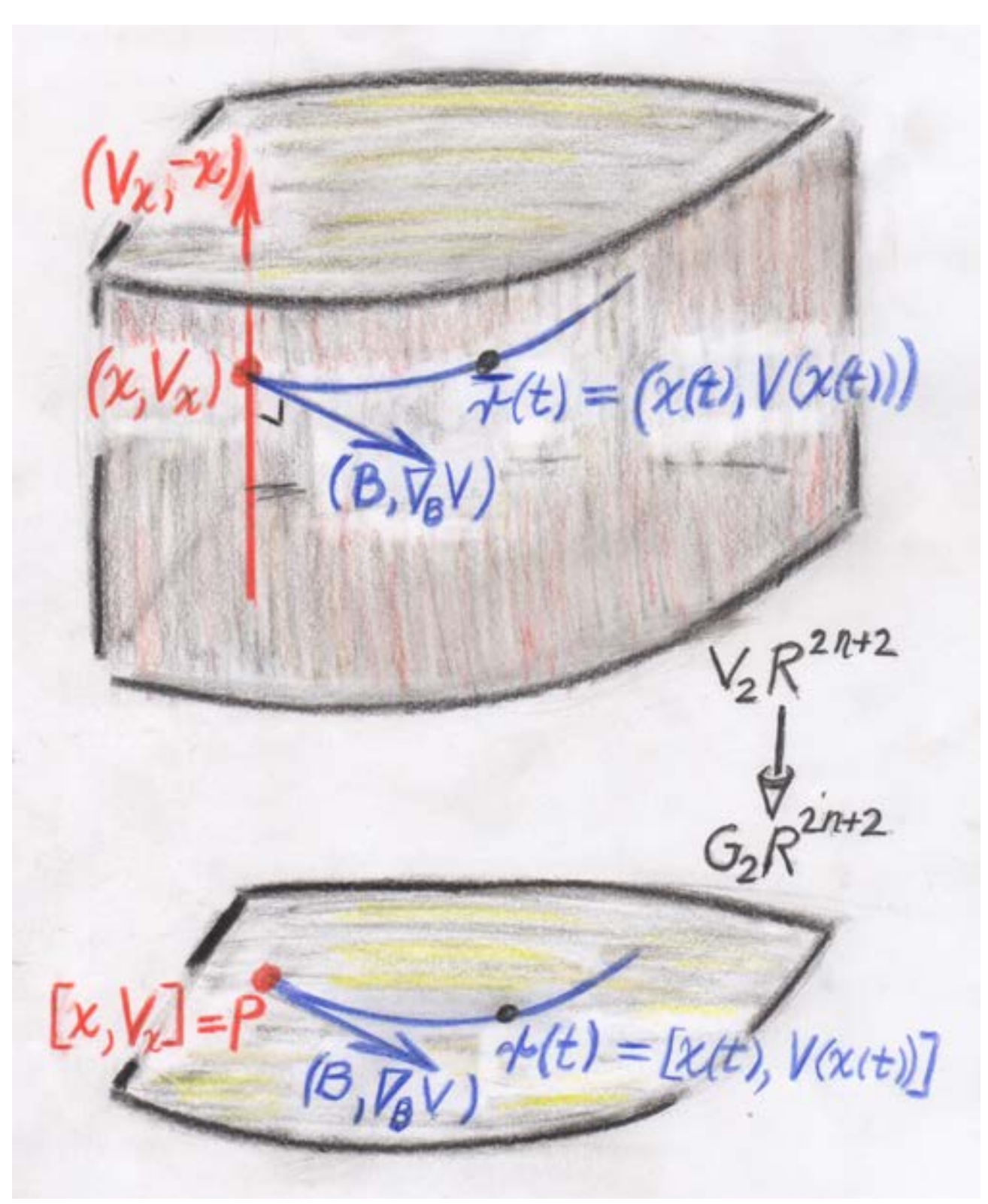

So this initial velocity $( B , \nabla_B V )$ of the curve $\bar{\gamma}(t) = ( x(t) , V(x(t)) )$ in the Stiefel manifold $V_2R^{2n+2}$ projects down undiminished to the initial velocity of the curve $\gamma(t) = [ x(t) , V(x(t)) ]$ in the Grassmann manifold $G_2R^{2n+2}$ .

The oriented 2-plane P with orthonormal basis x and $V_x$ begins to turn and twist so that x moves in the direction of B and so that $V_x$ moves in the direction of $\nabla_B V$ .

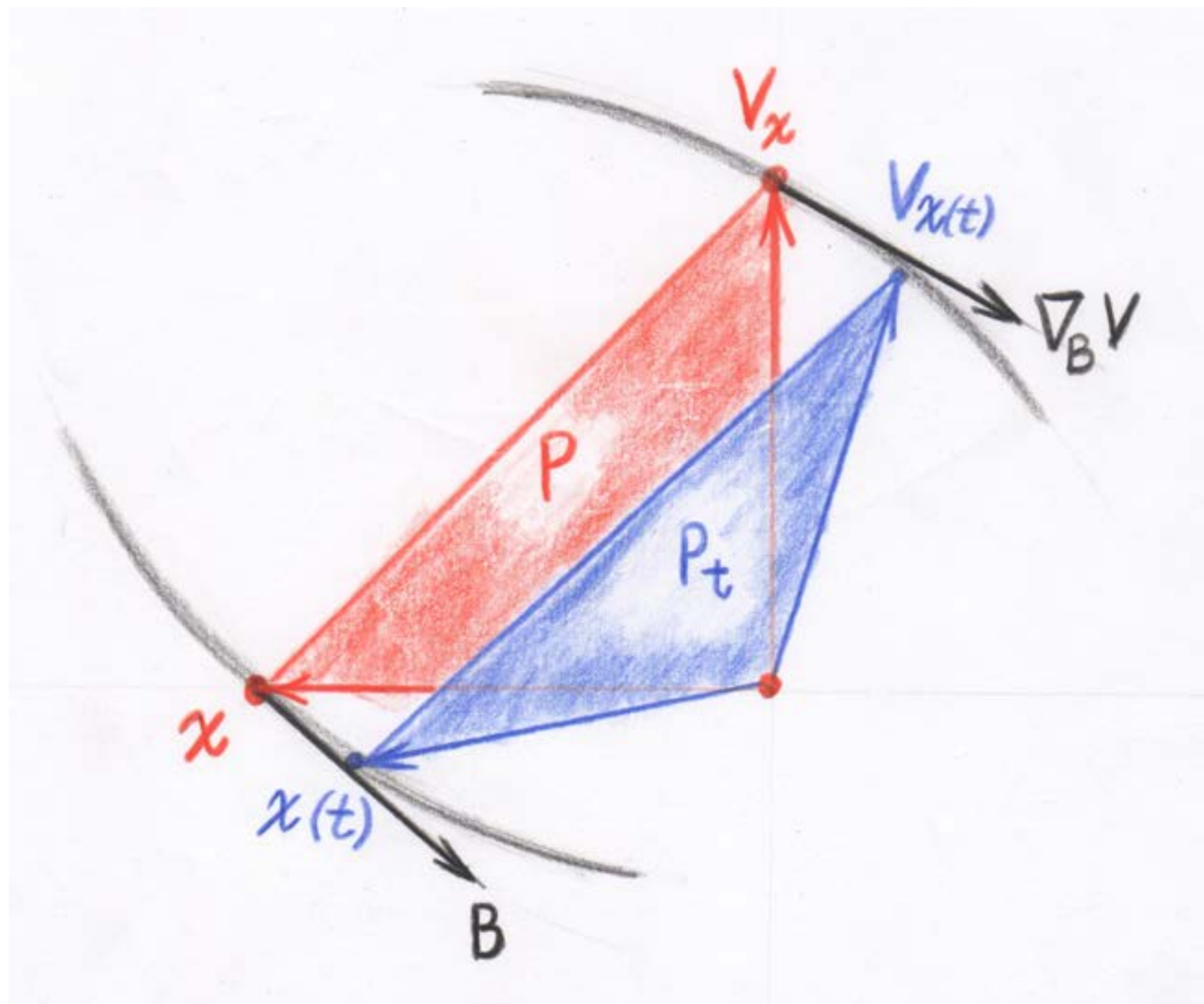

Since the curve $\gamma(t)$ runs within the base space $M_F$ of our fibration F of $S^{2n+1}$ by oriented great circles, we see that its initial velocity vector $( B , \nabla_B V )$ lies in the tangent space $T_P M_F$ to $M_F$ at P .

Since the vector B was an arbitrary choice of unit vector in $V_x^\perp$ , the vector $( B , \nabla_B V )$ will lie in $T_P M_F$ for all such choices of B in $V_x^\perp$ .

And then by linearity, the vectors $( B , \nabla_B V )$ will lie in $T_P M_F$ for all choices of B in $V_x^\perp$ , whether of unit length or not.

Counting dimensions, this must be all of $T_P M_F$ , which is therefore seen to be the graph of the twisting map $T_x = (\nabla V)|_{V_x^\perp} : V_x^\perp = P^\perp \rightarrow V_x^\perp = P^\perp$ , as claimed.

Hence by TOOL 2, the twisting map can have no real eigenvalues, completing the proof of our Proposition 1.

# CONTACT STRUCTURES

**PROPOSITION 2. Let $F : S^1 \subset S^{2n+1} -p\rightarrow M_F$ be a smooth fibration of the $2n+1$ sphere by oriented great circles, let $V$ be the unit vector field on $S^{2n+1}$ tangent to the fibres of $F$, and let $\xi_F = \{V_x^{\perp} : x \in S^{2n+1}\}$ be the distribution of tangent 2n-planes orthogonal to these fibres.**

**Let $T_x : V_x^{\perp} \rightarrow V_x^{\perp}$ be the twisting map introduced earlier.**

**Then a necessary and sufficient condition for $\xi_F$ to be a contact structure on $S^{2n+1}$ is that the skew-symmetrization $T_x - T_x^{tr}$ of this twisting map be non-singular at each point $x \in S^3$.**

**Proof.**

Let $\alpha$ be the differential one-form on $S^{2n+1}$ defined by $\alpha(W) = \langle V , W \rangle$ for each smooth vector field $W$ on $S^{2n+1}$.

By definition, $\xi_F$ will be a contact structure on $S^{2n+1}$ if and only if $\alpha \wedge (d\alpha)^n \neq 0$ at each point $x \in S^{2n+1}$, equivalently, if and only if the 2-form $d\alpha$ is non-degenerate on each 2n-plane $V_x^{\perp}$, meaning that for each nonzero vector $v \in V_x^{\perp}$, there is a vector $w \in V_x^{\perp}$ such that $d\alpha\,(v , w) \neq 0$.

The following argument, suggested by the referee, is a big improvement over our earlier coordinate-laden proof.

Let $X$ and $Y$ be vector fields on $S^{2n+1}$ which are orthogonal to $V$, or, in other words, lie in $\ker(\alpha)$.

Then

$$d\alpha(X, Y) = X(\alpha(Y)) - Y(\alpha(X)) - \alpha([X, Y]) \qquad \textit{Cartan Formula}$$

$$= X(\langle V, Y \rangle) - Y(\langle V, X \rangle) - \langle V, [X, Y] \rangle$$

$$= X(\langle V, Y \rangle) - Y(\langle V, X \rangle) - \langle V, \nabla_X Y - \nabla_Y X \rangle$$

*symmetric connection*

$$= \langle \nabla_X V, Y \rangle + \langle V, \nabla_X Y \rangle - \langle \nabla_Y V, X \rangle - \langle V, \nabla_Y X \rangle$$

$$- \langle V, \nabla_X Y \rangle + \langle V, \nabla_Y X \rangle$$

$$= \langle \nabla_X V, Y \rangle - \langle \nabla_Y V, X \rangle \qquad \textit{cancelling terms}$$

$$= \langle (\nabla V)(X), Y \rangle - \langle (\nabla V)(Y), X \rangle \qquad \textit{just rewriting}$$

$$= \langle X, (\nabla V)^{tr}(Y) \rangle - \langle X, (\nabla V)(Y) \rangle$$

$$= \langle X, ((\nabla V)^{tr} - \nabla V)(Y) \rangle$$

$$= \langle X, (T^{tr} - T)(Y) \rangle .$$

Copying the final simplification,

$$d\alpha(X, Y) = \langle X, (T^{tr} - T)(Y) \rangle .$$

If $T_x - T_x^{tr} : V_x^{\perp} \to V_x^{\perp}$ is nonsingular for each point $x \in S^{2n+1}$, then given any nonzero $X$ in $V_x^{\perp}$, we simply choose $Y$ in $V_x^{\perp}$ so that $(T^{tr} - T)(Y) = X$, and get

$$d\alpha(X, Y) = \langle X, (T^{tr} - T)(Y) \rangle = \langle X, X \rangle \neq 0 ,$$

confirming that $d\alpha$ is non-degenerate at $x$, and hence that $\xi_F$ is a contact structure on $S^{2n+1}$.

Contrariwise, if $T_x - T_x^{tr} : V_x^{\perp} \to V_x^{\perp}$ is singular for some $x \in S^{2n+1}$, we just choose a nonzero vector $Y$ in the kernel of $T_x - T_x^{tr}$, and get

$$d\alpha(X, Y) = \langle X, (T_x - T_x^{tr})(Y) \rangle = \langle X, 0 \rangle = 0$$

for all $X$ in $V_x^{\perp}$. This tells us that $d\alpha$ is degenerate at $x$, and hence that $\xi_F$ fails to be a contact structure on $S^{2n+1}$.

## EXAMPLES AND COUNTER-EXAMPLES

**The exceptional case on the 3-sphere: for every smooth great circle fibration, the tangent hyperplane distribution orthogonal to its fibres is a contact structure, which in fact is tight.**

For $n = 1$, we are looking at 2 x 2 matrices $T =$

$$\begin{matrix} a_{11} & a_{12} \\ a_{21} & a_{22} \end{matrix}$$

For $T$ to have no real eigenvalues, we must have

$$(Tr\, T)^2 - 4\, det\, T < 0 .$$

Writing this out,

$$(a_{11} + a_{22})^2 - 4\,(a_{11}\, a_{22} - a_{12}\, a_{21}) < 0 ,$$

equivalently,

$$(a_{11} - a_{22})^2 < -4\, a_{12}\, a_{21} .$$

The left side is $\geq 0$, and hence $a_{12}\, a_{21} < 0$. So one of the off diagonal terms must be positive and the other negative, and hence their difference $a_{12} - a_{21}$ cannot equal zero.

Thus $T - T^{tr} =$

$$\begin{matrix} 0 & a_{12} - a_{21} \\ a_{21} - a_{12} & 0 \end{matrix}$$

is non-singular.

It follows from Proposition 2 that for any great circle fibration of the 3-sphere, the orthogonal distribution of tangent 2-planes must be a contact structure.

**It remains to see why this contact structure must be tight.**

Following [Gluck-Warner 1983] and [Gluck 2022], let $F$ be a smooth great circle fibration of $S^3$ , and $H$ the Hopf fibration to which it is connected in the deformation retraction of the space of all great circle fibrations of $S^3$ to its subspace of Hopf fibrations.

This deformation retraction provides a one-parameter family $F_t$ of such fibrations, which begins with $F$ at $t = 0$ and ends with $H$ at $t = 1$ .

Then the corresponding contact forms $\alpha$ for $F$ and $\alpha'$ for $H$ can also be connected by a one-parameter family $\alpha_t$ of contact forms.

Hence by the Gray Stability Theorem [Gray, 1959; Geiges, 2008]], there is an isotopy $h_t$ of diffeomorphisms of $S^3$ with $h_0 = \text{identity}$ and with $h_t^*(\alpha_0) = f(t)\, \alpha_t$ , where $f(t)$ is a real-valued function.

Thus the contact structures $\xi_F$ and $\xi_H$ are isotopic, meaning that there is diffeomorphism $h\colon S^3 \to S^3$ , isotopic to the identity, such that $dh(\xi_F) = \xi_H$ .

Since $\xi_H$ is tight, so also is $\xi_F$ tight.

This reproves our earlier result without using special low-dimensional tools.

## A counterexample on the 5-sphere.

We construct now a smooth fibration of the 5-sphere by great circles whose orthogonal 4-plane distribution is not a contact structure.

For $n = 2$, we are looking at 4 x 4 matrices, and we give a specific example of such a matrix $T$ with no real eigenvalues, for which $T - T^{tr}$ is singular.

Let

$$T = \begin{pmatrix} 0 & 1/2 & 1 & 0 \\ -1/2 & 0 & 0 & 1 \\ 0 & 0 & 0 & 1/2 \\ 0 & 0 & -1/2 & 0 \end{pmatrix}.$$

The eigenvalues of $T$ are the purely imaginary numbers $i/2$ and $-i/2$, each of multiplicity two.

Then

$$T - T^{tr} = \begin{pmatrix} 0 & 1 & 1 & 0 \\ -1 & 0 & 0 & 1 \\ -1 & 0 & 0 & 1 \\ 0 & -1 & -1 & 0 \end{pmatrix}.$$

So we have on our hands a real 4 x 4 matrix $T$ with no real eigenvalues, for which $T - T^{tr}$ is singular.

From this preassigned twisting map we construct a germ of a fibration of the 5-sphere by oriented circles, thanks to TOOL 2 and the local version of TOOL 1, and then extend this germ to a fibration of the entire 5-sphere by oriented great circles, thanks to TOOL 3, and in that way obtain our counterexample.

## Counterexamples on the remaining odd-dimensional spheres.

Guided by the preceding section, we now construct a smooth fibration of the 2n+1 sphere by oriented great circles whose orthogonal 2n-plane distribution is not a contact structure.

We need a 2n x 2n matrix $T$ with no real eigenvalues, for which $T - T^{tr}$ is singular.

To get this, we place 2 x 2 blocks

$$\begin{matrix} 0 & 1/2 \\ -1/2 & 0 \end{matrix}$$

down the diagonal, using $n$ of them, and then a single 2 x 2 block

$$\begin{matrix} 1 & 0 \\ 0 & 1 \end{matrix}$$

in the upper right corner.

The eigenvalues of this matrix $T$ are the purely imaginary numbers $i/2$ and $-i/2$, this time each of multiplicity $n$.

Then $T - T^{tr}$ has 2 x 2 blocks

$$\begin{matrix} 0 & 1 \\ -1 & 0 \end{matrix}$$

down the diagonal, using $n$ of them, a single 2 x 2 block

$$\begin{matrix} 1 & 0 \\ 0 & 1 \end{matrix}$$

in the upper right corner, and a single 2 x 2 block

$$\begin{matrix} -1 & 0 \\ 0 & -1 \end{matrix}$$

in the lower left corner.

The first and last rows of $T - T^{tr}$ are negatives of one another, while the second and next-to-last rows of $T - T^{tr}$ are identical. So $T - T^{tr}$ is singular.

From this matrix $T$ we construct a germ of a fibration of the 2n+1 sphere by oriented great circles, using TOOL 2 and the local version of TOOL 1 as in the preceding section, and then extend this germ to a fibration of the entire 2n+1 sphere by oriented great circles, once again using TOOL 3, and we have our desired counterexample.

## REFERENCES

1931 Heinz Hopf, *Über die Abbildungen der dreidimensionalen Sphäre auf die Kugelfläche, Math. Ann. 104, 637-665.*

1935 Heinz Hopf, *Über die Abbildungen von Sphären auf Sphären niedrigerer Dimension,* Fund. Math. 25, 427-440.

1959 John W. Gray, *Some global properties of contact structures*, Annals of Math, 69, 421 - 450.

1983 Herman Gluck and Frank Warner, *Great circle fibrations of the three-sphere,* Duke Math. Journal 50 (1), 107-132.

1983 Herman Gluck, Frank Warner and C. T. Yang, *Division algebras, fibrations of spheres by great spheres and the topological determination of space by the gross behavior of its geodesics,* Duke Math. Journal, Vol. 50 No. 4, 1041-1076.

1992 Manfredo do Carmo, *Riemannian Geometry,* Birkhäuser, Boston-Basel-Berlin.

1992 Yakov Eliashberg, *Contact 3-manifolds 20 years since J. Martinet's work,* Annales de L'Institut Fourier, vol 42, No. 1-2, 165-192.

1993 Yakov Eliashberg, *Classification of contact structures on* $R^3$ *,* International Math. Research Notices, 3, 87-91.

2003 John Etnyre, *Introductory lectures on contact geometry,* Topology and Geometry of Manifolds, Proc. Sympos. Pure Math. 71, Amer. Math. Soc., 81-107.

2004 Benjamin McKay, *The Blaschke conjecture and great circle fibrations of spheres,* Amer. J. Math. Vol 126, No. 5, 1155-1191

2008 Hansgeorg Geiges, *An Introduction to Contact Topology,* Cambridge Univ. Press.

2015 Benjamin McKay, *A summary of progress on the Blaschke conjecture,* ICCM Not. 3, No. 2, 33-45.

2018 Patricia Cahn, Herman Gluck and Haggai Nuchi, *Deformation and Extension of Fibrations of Spheres by Great Circles,* Algebraic and Geometric Topology, vol. 18, No. 3, 1323-1360.

2022 Herman Gluck, *Great circle fibrations and contact structures on the 3-sphere,* Geometriae Dedicata 216:72.

University of Pennsylvania

*gluck@math.upenn.edu*
*jingyey@sas.upenn.edu*